# Optimization of frequency quantization

VN Tibabishev

Asvt51@narod.ru

**We obtain the functional defining the price and quality of sample readings of the generalized velocities. It is shown that the optimal sampling frequency, in the sense of minimizing the functional quality and price depends on the sampling of the upper cutoff frequency of the analog signal of the order of the generalized velocities measured by the generalized coordinates, the frequency properties of the analog input filter and a maximum sampling rate for analog-digital converter (ADC). An example of calculating the frequency quantization for two-tier ADC with an input RC filter.**

## 1. Introduction

The well-known sampling theorem specifies the sampling frequency of $F_1$ and is in many publications [1]. It is applicable only to send the generalized coordinate systems, for example, the deviation of the analog signal from the zero level in the absence of noise as the substitution frequencies. Sampling frequency is chosen from $F_1 > 2F_v$, where $F_v$ - cut-off frequency range of the analog signal. Sampling frequency of $F_1$ is called the minimum sampling rate of the zero-order generalized speed.

Let an almost periodic signal is the Fourier

$$y(t) = \sum_{k=0}^{k=\infty} A_k \cos(2\pi F_k t + \theta_k), -\infty < t < +\infty . \qquad (1)$$

Sampling theorem is proved only for signals that satisfy all harmonic amplitudes
$A_k \equiv 0$, if the frequency of $F_k > F_v$ .                                      (2)

It is known [2] that the simultaneous assignment of all the coordinates and velocities completely defines the state of one-dimensional and multidimensional systems. In this case it is assumed that enough to know information about the generalized coordinate and generalized velocity of the first derivative. If the generalized coordinates and generalized velocities are described by continuously differentiable functions, higher derivatives are the first derivative for any value of the current time. The situation changes if we are given discrete samples of the generalized coordinates at time intervals $\Delta t$. To find the approximate value of the first derivative on the interval $\Delta t$ requires at least one additional count in the middle of the interval. The minimum sampling frequency for estimating the derivative of the $k$ − th order associated with a minimum frequency of the zero order by the obvious relation $F_2 = 2(k + 1)F_v$ , where $k$ −order derivative, $0 \leq k \leq n$. It follows that knowledge of an approximate estimate of the first derivative is a necessary but not sufficient condition for an approximate estimate of higher order derivatives.

In this regard, the task of sampling the generalized coordinates and generalized velocities $k$ -th order one-dimensional case.

The quality (accuracy) estimates for the derivatives of the generalized coordinates with increasing frequency quantization, but the fee increases (price) as the volume of information being processed for the same duration of implementation. Necessary to determine the optimum sampling frequency $F_o(k) > F_2 > F_1$, where the price-quality rate reaches a minimum value for the evaluation of the derivative $k$ - th order.

In reality, the original signal with a limited range of $y(t)$ is given as a finite function. It is the product of the implementation model signal (1) defined on the whole line, and finite window function $O(t)$, which differs from zero only on the observation interval signal. The spectrum of the product of two functions $y(t) \cdot O(t)$ is the convolution of the spectra of the original signal and finite window function [3]. The convolution of the spectra of the product of these functions is defined on the whole infinite frequency axis. Therefore, the condition (2) does not hold in the real world well-known theorem of quantization.

Analog low-pass filter mounted in front of an analog-digital converter (ADC) to approximate the condition of quantization of a theorem (2). The real low-pass filter may have limited the rate of decrease of the amplitude-frequency response. The cutoff frequency $f_s > F_{05}$ low-pass filter is chosen for a given level of suppression of high-frequency part of the spectrum, where $F_{05}$ - filter bandwidth at half power. So instead of frequency quantization $F_1$ , $F_2$ or $F_o(k)$ in the ADC is sampling frequency $F_p$, which substitute for the suppression of frequencies chosen from the condition $F_p > 2f_s$ After the suppression substitution frequencies continued use of the frequency quantization $F_p$ resulting in redundancy of information. This shortcoming is eliminated by the digital output low pass filter in the form of a frequency divider with division ratio equal $K_d = F_p/F_o(k)$.

## 2. Determining the optimal frequency of quantization of the generalized velocity

The condition of the limited spectrum for finite processes is approximately satisfied by the input analog low-pass filter. With the limited frequency spectrum, for example, $F_v$ there is a boundary harmonic component $y_v(t)$ with the highest incidence $\omega_v = 2\pi F_v$ of non-zero amplitude $A_v$, for example, $y_v(t) = A_v \cos(\omega_v t + \theta_v)$ which defines the boundary harmonic component of the generalized coordinates. It is believed that the amplitude of harmonics with frequencies $\omega_{v+m} > \omega_v$  $m = 1,2,3,...$ are so small that their influence can be neglected. In determining the quality (accuracy) estimate of the approximate values derived from their exact values will use the mean-square numerical characteristics of the species

$$J_i^2 = \lim_{T \to \infty} 1/T \int_{t=0}^{t=T} f_i^2(t) dt. \tag{3}$$

If you put in this functional $f_i(t) = A_v \cos(\omega_v t + \theta_v)$, then the value of the functional will not depend on the initial phase. Therefore, the value of the initial phase $\theta_v$ we assume to be zero in order to simplify the expressions obtained.

Boundary harmonic component of the generalized velocity $k$ - th order $y^{(k)}(t)$ coordinate $y(t)$ is the time derivative of $y^{(k)}(t) = d^k y(t)/(dt)^k$. With digital signal processing instead of an infinitely small $dt$ quantity may be taken only a finite quantity $\Delta t$, the maximum value which is determined by the sampling theorem for the frequency of $F_p$, and the minimum value is limited to the maximum sampling frequency (speed) $F_s$ applied to the ADC. With these constraints, the value $\Delta t$ can be written as, where is not necessarily a positive integer which can take values in the range where $\Delta t = 1/2F_p N$ where $N$, which can take values in the interval $[1, N_m]$, where $N_m = F_s/2F_p$ - the maximum allowable value.

An approximate estimate of the harmonic component of the generalized velocity $v^{(k)}(t)$ $k$-th order is given by

$$v^{(k)}(t, \Delta t) = \Delta^{(k)} y(t, \Delta t)/(\Delta t)^k,$$

where $\Delta^{(k)} y(t, \Delta t)$ - the finite difference $k$ - th order continuous function $y(t)$ on the entire line.

Known connection with the finite difference derivatives of the corresponding order [4]

$$\Delta^{(k)} y(t) = \Delta t^k y^{(k)}(t, t + ak\Delta t), 0 < a < 1.$$

We choose α such that would $\alpha k = 1$. At the same time

$$\Delta t = 1/2F_p N = \pi/\omega_p N.$$

Deviation of the approximate estimate of $v^{(k)}(t,\Delta t)$ from its exact value of $y^{(k)}(t)$ is given by $R(\xi) = v^{(k)}(\xi) - y^{(k)}(t)$.

As an absolute measure of quality assessment of the boundary of the harmonic component of the generalized velocities using the functional form (3)

$$J_1^2(R^2(\xi)) = \lim_{T\to\infty} 1/T \int_{t=0}^{t=T} R^2(\xi)\, dt =$$

$$= \omega_p^{2k} A^2 \cos^2(\omega_p t)(1-\cos(\omega_p \Delta t))^2/2.$$

We define the functional significance for the exact derivative of $k$ - th order harmonic component of the boundary at a frequency $\omega_p$

$$J^2\left(y_{\omega_p}^{(2k)}(t)\right) = A^2 \omega_p^{2k}/2.$$

We find a numerical characteristic of the relative error of estimate of the boundary of the harmonic component of the derivative $k$ - th order

$$r(k, N) = \frac{R(\xi)}{J_y} = 1 - \cos(\omega_p \Delta t), \tag{4}$$

where $\Delta t = 1/2F_p N = \pi/\omega_p N$, $0 < k < n$, $1 < N < F_s/2F_p$.

With the increasing $N$ importance of decreasing the value of the functional $J_1(k,N)$, but it increases the number of counts per unit time, which acts as a fee (prices) to improve the quality of measurement. If at a frequency quantization sampling theorem the number of counts per unit time $G_1 = 2F_p$, with increasing values of the number of counts per unit time $N > 1$ increases to a maximum value

$$G_m = F_m = 2F_p N_m = F_s.$$

As a functional measure takes into account both the quality of the boundary of the generalized harmonic component of speed and cost (price) for the quality of measurements can be taken form the functional

$$J(N) = r(k,N) + J_2(N), \tag{5}$$

where $J_2(N)$ - the functional defining the cost (price) for measuring the quality of the boundary of the harmonic component of the generalized velocities.

Functional $r(k,N) \in [0,1]$ is the relative numerical characteristic. Therefore, the functional $J_2(N)$ must also be defined as the relative numerical characteristic $J_2(N) \in [0.1]$. As such a functional can be, for example, the functional form

$$J_2(N) = G_1/G_m = 2NF_p/F_s. \tag{6}$$

At the same time
$$J(N) = 1 - \cos(\pi/N) + 2NF_p/F_s. \tag{7}$$

Number $N_o$ for which the functional $J(N_o)$ attains its minimum value determines the optimal sampling frequency $F_o = 2N_o F_p$ for estimating the boundary of the harmonic component of the generalized velocities in terms of minimum error for the minimum sampling rate at which we can neglect the noise substitution frequencies.

Found functionality is a transcendental equation. This circumstance does not allow to obtain an expression for the parameter $N$ in an explicit form, in which the functional $J(N)$ attains its minimum value. In this regard, evaluation of optimal sampling frequency is a numerical method.

An example. Let the spectrum of the source signal frequency is limited $F_v = 2$ kHz. The task is to determine the optimum sampling frequency for estimating the first order derivative in the sense of minimizing the functional "price - quality" (7). The minimum sampling frequency is $F_2 = 2(k+1)2\text{кHz} = 8\text{кHz}$ for estimating the generalized first-order rate.

Consider the order of solving the problem for two-link input RC = T low-pass filter. For a frequency of 2 kHz determined by the filter time T on the half-power $T = 5{,}12 \cdot 10^{-5}$ s. Cutoff frequency bandwidth of the convolution of the spectra is determined by the amplitude $f_s$ frequency characteristic of two-link RC-low pass filter on the selected level, for example, 0.01 signal power. It is $f_s = 4{,}56$ kHz. Assume that the maximum frequency conversion in the ADC F_m (speed) is 500 kHz.

Numerical method is established that the minimum of functional (7) is reached for $N = 5$. Hence, the optimal sampling frequency $F_o = 5 * 2 * 4.56 = 45.6$ kHz. Standard error of change in the generalized first-order rate is 4.8%. The frequency of the output pulses for ADC evaluation of the generalized first-order rate is 2 * 4.56 = 9.12 kHz. Coefficient Decimation is 5.

When these data to estimate the generalized coordinates must be taken in the optimal sampling frequency of the ADC is equal to 36.55 kHz sampling frequency instead of sampling theorem 4 kHz. The coefficient of decimation is 8. This is a fee for the suppression of the substitution frequencies. If you need information about the generalized coordinate and generalized velocity of the first order, then received the two frequencies of quantization is chosen the largest value, and the two values of the coefficients of decimation selects the smallest value.

For more information on the quantization of the signals in the solution of problems of identification of the dynamics of multidimensional control objects can be found at http://asvt51.narod.ru/

## Conclusion

From the known sampling theorem sampling rate determined by the inequality $f_k > 2f_v$ under two conditions. First, it is estimated only a generalized coordinate, and, secondly, it is believed that there is no substitute for digitizing frequency. In many cases it is necessary to digitize the generalized velocities of finite orders and eliminate the phenomenon of substitution frequencies. It is established that there exists an optimal sampling frequency in the sense of minimizing the functional, which determines the price and quality measurement. This frequency depends on the upper cutoff frequency $f_v$, of the signal the order of the derivative (generalized velocity), the frequency properties of the input analog low-pass filter and the maximum frequency reference to the ADC. Is functional, which determines the price and quality measurements. An example is given for calculating the optimum frequency quantization for two-link RC-low pass filter.

15/10/2011　　　　　　　　　　　　　　　　V. Tibabishev